\documentclass{amsart}
\usepackage{amssymb}
\usepackage[square]{natbib}
\usepackage{epsfig}
\newtheorem{mainthm}{Theorem}[]
\newtheorem{thm}{Theorem}[section]

\newtheorem*{thm*}{Theorem} 
\newtheorem{cor}[thm]{Corollary}
\newtheorem{lem}[thm]{Lemma}

\newtheorem{prop}[thm]{Proposition}

\theoremstyle{definition}

\newtheorem*{rems*}{Remarks}
\theoremstyle{remark}
\newtheorem{rem}[thm]{Remark}

\newcommand{\N}{\mathbb{N}} 

\newcommand{\R}{\mathbb{R}}

\newcommand{\tg}{\tilde{g}}
\newcommand{\Sph}{\mathbb{S}}
 \newcommand{\bh}{\bar{h}}

\DeclareMathOperator{\tri}{tri}\DeclareMathOperator{\tr}{tr}

\DeclareMathOperator{\Ric}{Ric}\DeclareMathOperator{\ad}{ad}
\newcommand{\wed}{\wedge}\DeclareMathOperator{\li}{l}
\DeclareMathOperator{\scal}{scal}

\DeclareMathOperator{\Or}{O}
\DeclareMathOperator{\SO}{SO}

\DeclareMathOperator{\Rc}{R}\DeclareMathOperator{\Sc}{S}
\DeclareMathOperator{\Wc}{W} 

\newcommand{\so}{\mathfrak{so}}

\DeclareMathOperator{\trace}{tr}

\DeclareMathOperator{\id}{id} 
\DeclareMathOperator{\Weyl}{W}

 \DeclareMathOperator{\diam}{diam}

\newcommand{\ml}{\langle}                     
\newcommand{\mr}{\rangle}                     

\newcommand{\eps}{\varepsilon}

\newcommand{\blambda}{\bar{\lambda}}

\newcommand{\bg}{\bar{g}}

\begin{document}
\thanks{The first author was supported by the Deutsche
Forschungsgemeinschaft}

\begin{titlepage}\title
{Manifolds with positive curvature operators are space forms}
\author{Christoph B\"ohm}\author{Burkhard Wilking}
\end{titlepage}
\maketitle
\setcounter{page}{1}
\setcounter{tocdepth}{0}


The Ricci flow has been introduced by Hamilton in 1982 \cite{H1}
in order to prove that a compact three-manifold admitting a
Riemannian metric of positive Ricci curvature is a spherical
space form. In dimension four Hamilton showed that  compact
four-manifolds with positive curvature operators are
spherical space forms as well \cite{H2}. More generally,
the same conclusion holds for compact four-manifolds with
$2$-positive curvature operators \cite{Che}. Recall that a curvature 
operator is called $2$-positive, if the sum of its two smallest 
eigenvalues is positive. In arbitrary dimen\-sions Huisken \cite{Hu} 
described an explicit open cone
in the space of curvature operators
such that the normalized Ricci flow evolves metrics whose
curvature operators are contained in that cone into metrics
of constant positive sectional curvature.

Hamilton conjectured that in all dimensions
compact Riemannian manifolds with positive curvature operators must
be space forms. In this paper we confirm this conjecture. 
More generally, we show the following

\begin{mainthm}\label{thm1}\label{mainthm}\label{mainthm1}
On a compact manifold the normalized Ricci flow evolves a Riemannian metric
with $2$-positive curvature operator to
a limit metric with constant sectional curvature.
\end{mainthm}

The theorem is known in dimensions below five \cite{H4}, \cite{H1}, \cite{Che}.
Our proof works in dimensions above two:
we only use Hamilton's maximum principle
and Klingenberg's injectivity radius estimate for quarter pinched 
manifolds. 
Since in dimensions above two a quarter pinched orbifold is covered 
by a manifold (see Proposition~\ref{prop: orbifold}),
our proof carries over to orbifolds.
 
This is no longer true in dimension two. In the manifold case
it is known that the normalized Ricci flow converges to a metric of constant curvature for any initial metric \cite{H4}, \cite{Cho}. However, there 
exist two-dimensional orbifolds with positive
sectional curvature which are not covered by a manifold.
On such orbifolds the Ricci flow converges to a nontrivial Ricci 
soliton \cite{CW}.

Let us  mention that a $2$-positive curvature operator has 
positive isotropic curvature.  Micallef and Moore \cite{MM} showed that
a simply connected compact manifold with positive isotropic curvature 
is a homotopy sphere.  However,  their techniques do not allow to get 
restrictions for the fundamental groups or the differentiable structure of
the underlying manifold.

We turn to the proof of Theorem \ref{thm1}. The (unnormalized)
Ricci flow is the geometric evolution equation
\begin{eqnarray*}
  \frac{\partial g}{\partial t}= -2\Ric(g)
\end{eqnarray*}
for a curve $g_t$ of Riemannian metrics on a compact manifold $M^n$.
Using moving frames, this leads to the following
evolution equation for the curvature operator $\Rc_t$ of $g_t$
(cf.~\cite{H2}):
\begin{eqnarray*}
   \frac{\partial \Rc}{\partial t}= \Delta \Rc + 2(\Rc^2+\Rc^\#)\,.
\end{eqnarray*}
Here $\Rc_t\colon \Lambda^2T_pM\rightarrow \Lambda^2T_pM$
and identifying $\Lambda^2T_pM$ with $\so(T_pM)$ we have
$$
  \Rc^\#=\ad \circ  \,(\Rc \wed \Rc) \circ\ad^*\,,
$$ 
where $\ad\colon \Lambda^2(\so(T_pM))\rightarrow \so(T_pM)$ is the adjoint 
representation.
Notice that in our setting the curvature operator of the round 
sphere of radius one is the identity.

We denote by $S^2_B(\so(n))$ the vectorspace 
of curvature operators, that is the vectorspace
of selfadjoint endomorphisms of $\so(n)$ 
satisfying the Bianchi identity.
Hamilton's maximum principle asserts
that a closed convex $\Or(n)$-invariant  subset $C$ of $S^2_B(\so(n))$ 
which is invariant under the ordinary differential equation
\begin{eqnarray}
   \frac{d\Rc}{dt}=\Rc^2+ \Rc^\#\label{ode}
\end{eqnarray}
defines a Ricci flow invariant curvature condition; 
that is, the Ricci flow evolves metrics on compact manifolds 
whose curvature operators at each point are contained in $C$ 
into metrics with the same property. 

In dimensions above four there are relatively 
few applications of the maximum principle, since
 in these dimensions the
ordinary differential equation~\eqref{ode} is not well understood. 
 By analyzing how the differential equation changes under linear equivariant 
transformations,
 we provide a general method  for constructing 
new invariant curvature conditions from known ones.

Any equivariant linear transformation of the space of curvature operators
respects the decomposition 
$$
 S^2_B(\so(n))= \ml I\mr \oplus \ml \Ric_0 \mr \oplus \ml \Weyl
\mr 
$$
into pairwise inequivalent irreducible $\Or(n)$-invariant subspaces.
Here $\ml I\mr$ denotes
multiples of the identity, $\ml \Weyl \mr$ the space of Weyl curvature
operators and $\ml \Ric_0 \mr$ are
the curvature operators of traceless Ricci type.
Given a curvature operator $\Rc$ we let $\Rc_{I}$ and
$\Rc_{\Ric_0}$ denote the projections onto $\ml
I\mr$ and $\ml \Ric_0 \mr$, respectively.
Furthermore let $\Ric\colon \R^n\rightarrow \R^n$ denote the Ricci tensor of $\Rc$
and $\Ric_0$ the traceless part of $\Ric$.

\begin{mainthm}\label{mainthm2}
For $a,b\in \R$ consider the equivariant linear map
\begin{eqnarray*}
  \li_{a,b}\colon S^2_B(\so(n))\to S^2_B(\so(n))\,\,\,;\,\,\Rc \mapsto 
   \Rc +  2(n-1)a \Rc_I+  (n-2)b\Rc_{\Ric_0}
\end{eqnarray*}
and let
\begin{eqnarray*}
D_{a,b}&:=&
   \li_{a,b}^{-1}\bigl(  (\li_{a,b}\Rc)^2+
    (\li_{a,b}\Rc)^\#\bigr)-\Rc^2-\Rc^\#\,.
\end{eqnarray*}
Then
\begin{eqnarray*}
D_{a,b}
&=&\bigl((n-2)b^2-2(a-b)\bigr) \Ric_0\wed \Ric_0
   +2 a \Ric \wed \Ric+\,2b^2 \Ric_0^2 \wed \id\\
&&+\frac{\tr(\Ric_0^2)}{n+2n(n-1)a}\bigl(nb^2(1-2b)-2(a-b)(1-2b+nb^2)\bigr)\, I\,.
\end{eqnarray*}
\end{mainthm}

The key fact about the difference $D_{a,b}$ of 
the pulled back differential equation and the differential equation
itself is that it does not depend on the Weyl curvature.

Let us now explain why Theorem \ref{mainthm2} 
allows us to construct new curvature conditions which are
invariant under the ordinary differential equation (\ref{ode}): 
We consider the image 
of a known invariant curvature condition $C$ under the linear map
$\li_{a,b}$ for suitable constants $a,b$. This new
curvature condition is invariant under the ordinary differential 
equation, if $\li_{a,b}^{-1}\bigl(  (\li_{a,b}\Rc)^2+(\li_{a,b}\Rc)^\#\bigr)$
lies in the tangent cone $T_{\Rc} C$ 
of the known
invariant set $C$. By assumption $\Rc^2+\Rc^\#$ lies 
in that tangent cone, and hence it suffices to show 
$D_{a,b}\in T_RC$. Since
this difference does not depend on the Weyl curvature, 
it can be solely computed from the Ricci tensor.

Using this technique
we construct a continuous family of invariant cones
joining the invariant cone of $2$-positive curvature operators
and the invariant cone of positive 
multiples of the identity operator. Then a standard
ode-argument shows that from any such family a generalized
pinching set can be constructed -- a concept which is slightly more general 
than Hamiltons concept of pinching sets in [H2].
In Theorem \ref{thm: from ode to pde} we show that Hamilton's 
convergence result carries over to our situation, completing the proof of
Theorem
\ref{thm1}.

We expect that  Theorem \ref{mainthm2}  and its K\"ahler analogue should 
give rise to further applications. This will be the subject 
of a forthcoming paper.


\section{Algebraic Preliminaries}\label{sec: prelim}

For a Euclidean vector space $V$ we let $\Lambda^2V$ denote the
exterior product of $V$. We endow $\Lambda^2 V$ with its natural
scalar product; if $e_1,\ldots, e_n$ is an orthonormal basis of
$V$ then $e_1\wed e_2,...,e_{n-1}\wed e_n$ is an orthonormal basis
of $\Lambda^2V$. Notice that two linear endomorphisms $A,B$ of $V$
induce a linear map
\begin{eqnarray*}
  A\wed B \colon \Lambda^2V \rightarrow \Lambda^2V\,\,;\,\,\,
   v\wed w\mapsto \tfrac{1}{2}\bigl( A(v)\wed B(w)+ B(v)\wed A(w)\bigr)\,.
\end{eqnarray*}
We will identify $\Lambda^2\R^n$ with the Lie algebra $\so(n)$ by
 mapping the unit vector $e_i\wed e_j$ onto the linear map
$L(e_i\wed e_j)$ of rank two which is a rotation with angle
$\pi/2$ in the plane spanned by $e_i$ and $e_j$. Notice that under
this identification the scalar product on $\so(n)$
corresponds to $\ml A,B\mr=-1/2\tr (AB)$.

For $n\geq 4$ there is a natural decomposition of
\begin{eqnarray*}
  S^2(\so(n))= \ml I \mr \oplus \ml \Ric_0 \mr \oplus \ml \Weyl \mr
               \oplus \Lambda^4(\R^n)
\end{eqnarray*}
into $\Or(n)$-invariant, irreducible and pairwise inequivalent 
subspaces.
An endomorphism $\Rc\in S^2(\so(n))$ satisfies the first Bianchi
identity if and only if $\Rc$ is an element in
$S^2_B(\so(n))= \ml I\mr \oplus \ml \Ric_0 \mr \oplus \ml \Weyl
\mr \,$. 
Given a curvature operator $\Rc\in S^2_B(\so(n))$ we let $\Rc_{I}$,
$\Rc_{\Ric_0}$ and $\Rc_{\Weyl}$, denote the projections onto $\ml
I\mr$, $\ml \Ric_0 \mr$ and $\ml \Weyl \mr $, respectively.
Moreover, let
$$
  \Ric\colon \R^n \to \R^n
$$
denote the Ricci tensor of $\Rc$, $\Ric_0$
the traceless Ricci tensor and
\begin{eqnarray}
   \blambda:= \tr (\Ric)/n   \quad \textrm{ and }\quad
   \sigma:= \|\Ric_0\|^2/n\,.\label{blambda}
\end{eqnarray}
Then
\begin{eqnarray}
    \Rc_I=\frac{\blambda}{n-1}\id \wed \id
    \quad \textrm{ and }\quad
     \Rc_{\Ric_0}=\frac{2}{n-2}\Ric_0 \wed \id\,.\label{RIRic}
\end{eqnarray}
Hamilton observed in \cite{H2} that next to the map
$(\Rc,\Sc)\mapsto \frac{1}{2}(\Rc\Sc+\Sc\Rc)$ there is a second
natural $\Or(n)$-equivariant bilinear map
$$
 \#\colon S^2(\so(n))\times S^2(\so(n))\to S^2(\so(n))
  \,\,;\,\,\,(\Rc,\Sc)\mapsto \Rc \# \Sc
$$
given by
\begin{eqnarray}
\langle (\Rc\#\Sc)(h),h\rangle
  &=&\frac{1}{2}
     \sum_{\alpha,\beta=1}^N \,\,\ml [\Rc(b_\alpha),\Sc(b_\beta)],h\mr
            \cdot    \ml [b_\alpha,b_\beta],h\mr\label{RSsh}
\end{eqnarray}
for $h\in \so(n)$ and an orthonormal basis $b_1,...,b_N$ of
$\so(n)$. The factor $1/2$ stems from that fact that we are using
the scalar product $-1/2\tr(AB)$ instead of $-\tr(AB)$ as in
\cite{H2}. We would like to mention that $\Rc \#\Sc =\Sc\#\Rc$ 
can be described invariantly
\begin{eqnarray*}
\Rc\#S= \ad \circ \,\,(\Rc\wed \Sc) \circ \ad^{*}\,,
\end{eqnarray*}
where $\ad\colon \Lambda^2\so(n)\rightarrow \so(n), u\wedge v
\mapsto [u,v]$ denotes the adjoint representation and $\ad^*$ is
its dual.
 Following Hamilton we set
\begin{eqnarray*}
  \Rc^\# =\Rc\#\Rc\,.
\end{eqnarray*}

We will also consider the trilinear form
\begin{eqnarray}\label{tri}
  \tri(\Rc_1,\Rc_2,\Rc_3)
    &=&\tr{ \big((\Rc_1 \Rc_2+ \Rc_2\Rc_1+2 \Rc_1\#\Rc_2)\cdot \Rc_3\big) }\,.
\end{eqnarray}
The authors learned from Huisken that $\tri$ is symmetric in all three
components.
 In fact using \eqref{RSsh} it is straightforward to check that
\[
\trace( 2(\Rc_1\#\Rc_2)\cdot \Rc_3)=\sum_{\alpha,\beta,\gamma=1}^N
\ml[\Rc_1(b_\alpha), \Rc_2(b_{\beta})], \Rc_3 (b_{\gamma})\mr\cdot
\ml[b_\alpha, b_{\beta}], b_{\gamma}\mr\,.
\]
Since the right hand side is clearly symmetric in all three
components this gives the desired result. Huisken also observed
that the ordinary differential equation \eqref{ode} is the
gradient flow of the function 
$$
  P(R)=\frac{1}{3}\tr (\Rc^3+\Rc \Rc^\#)=\frac{1}{6}\tri(\Rc,\Rc,\Rc)\,.
$$

Finally we
recall that if $e_1,\ldots,e_n$ denotes an orthonormal basis of eigenvectors
of $\Ric$, then
\begin{eqnarray}
\Ric(\Rc^2+\Rc^\#)_{ij}
   = \sum_{k}\Ric_{kk}\Rc_{kijk}\label{RicstrRicci}
\end{eqnarray}
where $\Rc_{kijk}=\ml\Rc(e_i\wed e_k),e_j\wed e_k\mr$, see \cite{H1}, \cite{H2}.


\section{A new Algebraic Identity for Curvature Operators}\label{sec: algid}

The main aim of this section is to prove Theorem~\ref{mainthm2}.
 A computation using
(\ref{RIRic}) shows that the linear map 
$\li_{a,b}\colon S^2_B(\so(n))\rightarrow S^2_B(\so(n))$  given
in Theorem~~\ref{mainthm2}
satisfies
\begin{eqnarray}\nonumber
\li_{a,b}(\Rc) 
&=&\Rc+ 2b\Ric\wed
   \id+2(n-1)(a-b)\Rc_I\label{li}\,.
\end{eqnarray}
The bilinear map $\#$ induces a linear $\Or(n)$-equivariant map given by
$\Rc\mapsto \Rc \#I$. The normalization 
of our parameters is related to the eigenvalues of this map.

\begin{lem}\label{lem: MshI}\label{lem: sharp eigenvalues}
Let $\Rc\in S^2_B(\so(n))$. Then
\begin{eqnarray*}
 \Rc + \Rc \#I = (n-1)\Rc_I+\frac{n-2}{2}\Rc_{\Ric_0}=\Ric \wedge \id\,.
\end{eqnarray*}
\end{lem}

\begin{proof}
One can write 
\begin{eqnarray}\label{polarization}
 \Rc+ \Rc \#I&=& 
   \tfrac{1}{4}\bigl( (\Rc+I)^2+(\Rc+I)^\#-(\Rc-I)^2-(\Rc-I)^\#\bigr)\,.
\end{eqnarray} 
The result on the eigenvalues of the map
corresponding to the subspaces $\ml \Ric_0 \mr$ and $\ml I\mr$
now follows from  equation \eqref{RicstrRicci}
 by a straightforward computation.
For $n=4$ one verifies directly that $\ml \Wc \mr$ is in the
kernel of the map  $\Rc\mapsto \Rc+ \Rc \#I$. Since there is a
natural embedding of the Weyl curvature operators in
$S^2_B(\so(4))$ to the Weyl curvature operators in $S^2_B(\so(n))$
this implies the same result for $n\ge 5$.
\end{proof}

We say that a curvature operator $\Rc$ is of Ricci type,
if $\Rc=\Rc_I+\Rc_{\Ric_0}$.

\begin{lem}\label{lem: RRiccitype}\label{lem: Riccitype}
Let $\Rc \in S^2_B(\so(n))$ be a curvature operator of Ricci type,
and let $\blambda$ and $\sigma$
be as in \eqref{blambda}. Then
\begin{eqnarray*}
  \Rc^2+\Rc^\#&=& \frac{1}{n-2}\Ric_0\wed \Ric_0
              +\frac{2\blambda}{(n-1)}\Ric_0\wed \id
         -\frac{2}{(n-2)^2}(\Ric_0^2)_0\wed \id\\
            &&   +\frac{\blambda^2}{n-1}I
           +\frac{\sigma}{n-2} I\,.
\end{eqnarray*}
Moreover
\begin{eqnarray*}
  \bigl(\Rc^2+\Rc^\#\bigr)_{\Wc}
   &=&
     \frac{1}{n-2}\bigl(\Ric_0 \wed \Ric_0\bigr)_{\Wc}\\
  \Ric( \Rc^2+\Rc^\#)
   &=&
     -\frac{2}{n-2} (\Ric_0^2)_0+\frac{n-2}{n-1}\blambda \Ric_0
     + \blambda^2 \id +\sigma \id\,.
 \end{eqnarray*} 
\end{lem}

\begin{proof} By equation~\eqref{RIRic}
$$
      \Rc= \Rc_I+ \Rc_{\Ric_0} 
         = \frac{\blambda}{(n-1)}I+\frac{2}{(n-2)} \Ric_0\wed \id\,.
$$
Using the abbreviation $\Rc_0=\Rc_{\Ric_0}$ we have 
\[
\Rc^2+\Rc^\#= \Rc^2_0+\Rc^\#_0+ \frac{2\blambda}{(n-1)}(\Rc_0+\Rc_0\#I)
+ \frac{\blambda^2}{(n-1)^2}(I+I^\#)\,.
\] 
Since the last two summands are known by Lemma~\ref{lem: sharp eigenvalues}, 
we may assume that 
$\Rc=\Rc_{\Ric_0}$.
Let $\lambda_1,\ldots,\lambda_n$ denote the eigenvalues of
$\Ric_0$ corresponding to an orthonormal basis $e_1,\ldots,e_n$ of
$\R^n$.
The curvature operator $\Rc$ is diagonal with respect to $e_1\wed e_2,...,e_{n-1}\wed e_n$ and we denote by $\Rc_{ij}=\tfrac{\lambda_i+\lambda_j}{n-2}$
 the corresponding eigenvalues for $1\leq i<j\leq n$.
Inspection of (\ref{RSsh}) shows that
also $\Rc^2+\Rc^\#$ is diagonal with respect to this basis. We have
\begin{eqnarray*} 
   (\Rc^2+\Rc^\#)_{ij}
      &=&\Rc_{ij}^2+\sum_{k\neq i,j}\Rc_{ik}\Rc_{jk}\\
      &=&\frac{(\lambda_i+\lambda_j)^2}{(n-2)^2} 
          +\frac{1}{(n-2)^2}
          \sum_{k\neq i,j}(\lambda_i+\lambda_k)(\lambda_j+\lambda_k)\\
      &=& \frac{\lambda_i\lambda_j}{(n-2)}
          +\frac{n\sigma-\lambda_i^2-\lambda_j^2}{(n-2)^2}\\
\end{eqnarray*}
as claimed.

The second identity follows immediately from the first.
To show the last identity notice that the Ricci tensor of 
$\Ric_0\wed \Ric_0$ is given by $-\Ric_0^2$. A computation shows
the claim.
\end{proof}

\begin{proof}[Proof of Theorem~\ref{mainthm2}.]
We first verify that $D=D_{a,b}$ does not depend on the Weyl
curvature of $\Rc$. We view $D$ as quadratic form in $\Rc$. Then
$$
  B(\Rc,\Sc):=\tfrac{1}{4} \bigr(D(\Rc+\Sc)-D(\Rc-\Sc)\bigl)
$$
is the corresponding bilinear form.

Let $\Sc=\Wc\in \ml \Weyl \mr$. We have to show $B(\Rc,\Wc)=0$ for all
$\Rc\in S^2_B(\so(n))$. We start by considering $\Rc\in \ml \Weyl
\mr$. Then $\li_{a,b}(\Rc\pm \Wc)=\Rc\pm \Wc$. It follows from formula
(\ref{RicstrRicci}) for the Ricci curvature of $\Rc^2+\Rc^\#$ that
$(\Rc\pm \Wc)^2+(\Rc\pm \Wc)^\#$ has vanishing Ricci tensor. Hence
$(\Rc\pm \Wc)^2+ (\Rc\pm \Wc)^\#$ is a Weyl curvature operator and
accordingly fixed by $\li_{a,b}^{-1}$.

Next we consider the case that $\Rc=I$ is the identity. Using the polarization
formula (\ref{polarization}) for $\Wc$ we see that
$ B(I,\Wc)$ is a multiple of $\Wc+\Wc\#I$, which is zero
by Lemma \ref{lem: MshI}.

It remains to consider the case of $\Rc\in \ml \Ric_0\mr$.
Using the symmetry of the trilinear form $\tri$ defined in
\eqref{tri}
we see for each $\Wc_2\in \ml \Wc\mr$ that
$$
\tri(\Wc,\Rc,\Wc_2)=\tri(\Wc,\Wc_2,\Rc)=0
$$
as $\Wc\Wc_2+\Wc_2 \Wc+2\Wc\#\Wc_2$ lies in $\ml \Weyl \mr$ and
$\Rc\in \ml \Ric_0\mr$. Combining this with $\tri(\Wc,\Rc,I)=0$ gives that
$\Wc\Rc+\Rc\Wc+2\Wc \#\Rc \in \ml \Ric_0\mr$. Using once
more that $\li:=\li_{a,b}$ is the identity on $\ml \Weyl \mr$ we
see that
$$
\li(\Wc)\li(\Rc)+\li(\Rc)\li(\Wc)+2\li(\Wc)\#\li(\Rc)
 =\li(\Wc\Rc+\Rc\Wc+2\Wc \#\Rc)\,.
$$
This clearly proves $B(\Rc,\Wc)=0$.

Thus, for computing $D$ we may assume that
$\Rc_{\Wc}=0$. So let $\Rc=\Rc_I+\Rc_{\Ric_0}$. 
We next verify that both sides of the equation 
have the same projection to the space $\ml \Wc \mr$ 
of Weyl curvature operators. 
Recall that $\li_{a,b}^{-1}$ induces the identity on 
$\ml \Wc\mr $ and that $\Ric_0(\li_{a,b}(R))=(1+(n-2)b)\Ric_0$.
Then using the second identity in Lemma~\ref{lem: Riccitype} we see that
\begin{eqnarray*}
D_{\Wc}
   &=& \frac{1}{n-2}( (1+(n-2)b)^2-1)
       \bigl(\Ric_0 \wed \Ric_0\bigr)_{W}\\
   &=& \bigl((n-2)b^2 +2b\bigr)
        \bigl(\Ric_0 \wed \Ric_0\bigr)_{W}\,.
\end{eqnarray*}
It is straightforward to check that the right hand side 
in the asserted identity for $D$ has the same projection 
to $\ml \Wc \mr$. 

It remains to check that both sides of the equation 
have the same Ricci tensor. 
Because of $\Ric(\li_{a,b}(\Rc))= (1+(n-2)b)\Ric_0+(1+2(n-1)a)\blambda \id$,
the third identity in Lemma~\ref{lem: Riccitype} implies
\begin{eqnarray}
\Ric(D)&=&-2b(\Ric_0^2)_0+2(n-2)a\blambda \Ric_0
+ 2(n-1)a\blambda^2 \id \nonumber\\[1ex]
&&+\frac{2(n-2)b+(n-2)^2b^2-2(n-1)a}{1+2(n-1)a}\nonumber\sigma \id\\[1ex]
&=&-2b\Ric_0^2+2(n-2)a\blambda \Ric_0
+ 2(n-1)a\blambda^2 \id \label{ricciformula}\\[1ex]
&&+\frac{2(n-1)b+(n-2)^2b^2-2(n-1)a(1-2b)}{1+2(n-1)a}\sigma \id\,.\nonumber
\end{eqnarray}
A straightforward computation shows that the same holds 
for the Ricci tensor of the right hand side 
in the asserted identity for $D$. This completes the proof. 
\end{proof}

\begin{cor}\label{cor: Deigen} We keep the notation of
Theorem~{\rm \ref{mainthm2}}, and let $\sigma$, $\blambda$ be as in \eqref{blambda}. Suppose that $\lambda_1,\ldots,\lambda_n$ are the eigenvalues of
$\Ric_0$ corresponding to an orthonormal basis $e_1,\ldots,e_n$.
Then $e_i\wed e_j$ $(i< j)$ is an eigenvector of
$D_{a,b}$ corresponding to the eigenvalue
\begin{eqnarray*}
d_{ij}&=&\bigl((n-2)b^2-2(a-b)\bigr) \lambda_i\lambda_j
          +2a(\blambda +\lambda_i)( \blambda +\lambda_j)
          +b^2(\lambda_i^2+\lambda_j^2)\\
    &&+\frac{\sigma}{1+2(n-1)a}\bigl(
       nb^2(1-2b)-2(a-b)(1-2b+nb^2)\bigr)\,.
\end{eqnarray*}
Furthermore,  $e_i$ is an eigenvector of the Ricci tensor of $D_{a,b}$
with respect to the eigenvalue
\begin{eqnarray*}
r_i
 &=&-2b\lambda_i^2
    +2a\blambda(n-2) \lambda_i+2a(n-1)\blambda^2\\
 & &+\frac{\sigma}{1+2(n-1)a}\bigl(n^2b^2-2(n-1)(a-b)(1-2b)\bigr)\,.
\end{eqnarray*}
\end{cor}

Notice that $\lambda_i+\blambda$ are the eigenvalues of
the Ricci tensor $\Ric$.
The first formula follows immediately from Theorem~\ref{mainthm2}, the second from (\ref{ricciformula}).


\section{New Invariant Sets}\label{sec: transversal}

We call a continuous family $C(s)_{s\in [0,1)}\subset
S_B^2(\so(n))$  of  closed convex $\Or(n)$-invariant cones of full
dimension a pinching family, if
\begin{enumerate}
\item each $\Rc\in C(s)\setminus \{0\}$ has positive scalar curvature,

\item $\Rc^2+\Rc^\#$ is contained in the interior of the tangent
      cone of $C(s)$ at $\Rc$ for all $\Rc \in C(s)\setminus \{0\}$
      and all $s\in (0,1)$,
\item $C(s)$ converges in the pointed Hausdorff topology to the
      one-dimensional cone $\R^+ I$ as $s\to 1$.
\end{enumerate}

The main aim of this section is to prove 

\begin{thm}\label{thm: pinching fam}
 There is a pinching family $C(s)_{s\in [0,1)}$ 
of closed convex cones
 such that $C(0)$ is the cone of $2$-nonnegative curvature operators.
\end{thm}

As before
 a curvature operator is called $2$-nonnegative if the sum
of its smallest two eigenvalues is nonnegative. 
It is known that the cone of $2$-nonnegative curvature 
operators is invariant under  the ordinary differential equation (\ref{ode})
(see \cite{H3}). 
The pinching family that we construct 
for this cone is defined piecewise by three subfamilies. 
Each cone in the first subfamily is the 
image of the cone of $2$-nonnegative curvature operators under
a linear map. In fact we have the following general result.

\begin{prop}\label{prop: newinv}
Let $C\subset S^2_B(\so(n))$ be a closed convex $\Or(n)$-invariant
subset which is invariant under the ordinary differential equation
{\rm (\ref{ode})}. Suppose that $C\setminus \{0\}$ is contained in
the half space of curvature operators with positive scalar
curvature, that each $\Rc\in C$ has nonnegative Ricci curvature
and that $C$ contains all nonnegative curvature
operators of rank 1. Then for $n\geq 3$ and
\begin{eqnarray*}
  b\in \big(0,\tfrac{\sqrt{2n(n-2)+4}-2}{n(n-2)}\big] \quad \textrm{ and }\quad
 2a=2b+(n-2)b^2
\end{eqnarray*}
the set $\li_{a,b}(C)$ is invariant under the vector field
corresponding to {\rm (\ref{ode})} as well. In fact, it
 is transverse to the boundary
of the set at all boundary points $\Rc\neq 0$.
\end{prop}

Using the Bianchi identity it is straightforward to check
that
a nonnegative curvature operator of rank 1
corresponds up to a positive factor and a change of basis in $\R^n$
to the curvature operator of $\Sph^2\times \R^{n-2}$.
The condition that $C$ contains all these
operators is equivalent to saying that $C$ contains
the cone of geometrically nonnegative curvature operators.
A curvature operator is geometrically nonnegative if it can be written as
the sum of nonnegative curvature operators of rank 1.
In dimensions above $4$ this cone is strictly smaller than the cone
of nonnegative curvature operators.
Although we will not need it, we remark that
the cone of geometrically nonnegative curvature operators
is invariant under {\rm (\ref{ode})} as well.

\begin{proof}
We have to prove
that for each $\Rc\in C$ the curvature operator
\begin{eqnarray}
   X_{a,b}&=&\li^{-1}_{a,b}(\li_{a,b}(\Rc)^2+\li_{a,b}(\Rc)^\#)\label{pullback}
\end{eqnarray}
lies in the tangent cone $T_{\Rc}C$ of $C$ at the point $\Rc$.
Notice that by assumption we have $\Rc^2+\Rc^\#\in T_{\Rc}C$. Thus
it suffices to show that $ D_{a,b}=X_{a,b}-\Rc^2-\Rc^\#$ lies in
$T_{\Rc}C$. Since $C$ contains all nonnegative
curvature operators of rank 1, we can establish this by showing
that $D_{a,b}$ is positive for $b>0$. Looking at the
formula for the eigenvalues of $D_{a,b}$ in Corollary \ref{cor: Deigen}
this amounts to showing that
\begin{eqnarray*}
0 &\le & b^2\big(n(1-2b)- (n-2)(1-2b+nb^2)\big)
\end{eqnarray*}
holds in the given range.
This is a straightforward computation.
\end{proof}

Let us remark that 
the intersection of two closed convex $\Or(n)$-invariant cones, which are
invariant under the ordinary differential equation (\ref{ode}), have the 
same properties as the given cones.

\begin{cor}\label{cor: 2positive} In order to prove  
Theorem~{\rm \ref{thm: pinching fam}}, it suffices 
to establish 
the existence of a pinching family $C(s)_{s\in [0,1)}$ 
with $C(0)$ being the cone of nonnegative curvature operators. 
\end{cor}

\begin{proof} Suppose $n\geq 4$. 
Notice that the cone $C$  of $2$-nonnegative curvature operators 
satisfies the assumptions of  Proposition~\ref{prop: newinv}. 
We plan to show that the family of closed invariant
cones  from Proposition~\ref{prop: newinv} can be extended
 to a pinching family. 
 By the above remark
it suffices to show that $\li_{a(b),b}(C\setminus\{0\})$ 
is contained in the open cone of positive curvature operators where 
$b$ is the maximal allowed value
from Proposition \ref{prop: newinv}.  In fact then we 
can extend the family from Proposition~\ref{prop: newinv} 
to a pinching family by defining it on the second part 
of the interval as a reparameterization of the pinching family 
 $(C(s)\cap C')_{s\in [0,1)}$ where $C':=\li_{a(b),b}(C)$.

Let $\Rc\in C\setminus\{0\}$. 
 Recall that by (\ref{li}) we have
$\li_{a,b}(\Rc)=\Rc+ 2b\Ric\wed \id+h\Rc_I$ for
$h:=2(n-1)(a-b)$.
The smallest eigenvalue of $\Rc$ is by a standard estimate larger
than or equal to $-\tfrac{2\tr(\Rc)}{n(n-1)-2}$. Moreover,
since the sum of the two smallest eigenvalues of $\Rc$ is 
nonnegative the smallest eigenvalue of $\Ric$ is bounded
from below by $(n-3)$ times the absolute value of the smallest eigenvalue of
$\Rc$. Thus in order to show that $\li_{a,b}(\Rc)>0$ it is sufficient to prove
$h>(1-2b)\tfrac{n(n-1)}{n(n-1)-2}$.
This is equivalent to
$$
  (n-2)b^2>(1-2b)\frac{n}{(n+1)(n-2)}\,.
$$
By the definition of $b$ we have $(n-2)b^2=\tfrac{2}{n}(1-2b)$.
This shows the claim for $n\geq 4$.  For $n=3$ 
Theorem~\ref{thm: pinching fam} is well known.
\end{proof}

It remains to construct a pinching family for the 
cone of nonnegative curvature operators. 
This pinching family 
will be defined up to parameterization 
piecewise by two subfamilies in the next two lemmas.

\begin{lem} \label{lem: firstdef}\label{lem: first family}
For $b\in [0,1/2]$ put
\begin{eqnarray*}
  a=\frac{(n-2)b^2+2b}{2+2(n-2)b^2}\quad \textrm{ and }\quad
  p=\frac{(n-2)b^2}{1+(n-2)b^2}\,.
\end{eqnarray*}
Then the set
$$
\li_{a,b}\Bigl(\bigl\{ \Rc\in S^2_B(\so(n))\bigm| \Rc\ge 0,
\Ric\ge p(b)\tfrac{\tr(\Ric)}{n}\bigr\}\Bigr)
$$
is invariant under the vector field corresponding to {\rm
(\ref{ode})}. In fact, for $b\in (0,1/2]$ it is transverse to the
boundary of the set at all boundary points $\Rc\neq 0$.
\end{lem}

\begin{proof}
Put
$$
 C(p):=\bigl\{ \Rc\in S^2_B(\so(n))\bigm| \Rc\ge 0, \Ric\ge
         p(b)\tfrac{\tr(\Ric)}{n}\bigr\}\,.
$$
It suffices to check that for $\Rc\in C(p)\setminus \{0\}$ the
pulled back vector field $X_{a,b}$ defined in (\ref{pullback})
is in the interior of the tangent cone of $C(p)$ at $\Rc$.

In the first step we verify that $X_{a,b}$ is positive definite
for $b\in (0,1/2]$. Since $\Rc^2+\Rc^\#$ is positive
semi-definite, we can establish $X_{a,b}> 0$  by showing $D_{a,b}>
0$.  Since by assumption $\Rc\in C(p)$ we have
the following estimate for the eigenvalues of $\Ric_0$:
$$
  \lambda_i \ge-(1-p)\blambda\,.
$$
Next, observe that
$$
  2(a-b)=\frac{1-2b}{1+(n-2)b^2}(n-2)b^2\,.
$$
We use the notation of Theorem~\ref{mainthm2} and
Corollary~\ref{cor: Deigen}. Rewriting $d_{ij}$ gives
\begin{eqnarray}
d_{ij}
  &=&\tfrac{2a}{1-p}( (1-p)\blambda +\lambda_i)
     ( (1-p)\blambda +\lambda_j)+2a p\blambda^2
      +b^2(\lambda_i^2+\lambda_j^2)\label{equ}\\
  & &+\frac{n(1+(n-2)b^2)
      -(n-2)(1-2b+nb^2)}{(1+2(n-1)a)(1+(n-2)b^2)}\sigma b^2(1-2b)\nonumber\\
   &>& \frac{2+2(n-2)b}{(1+2(n-1)a)(1+(n-2)b^2)}\sigma b^2(1-2b)
        \nonumber\\
  &\ge& 0\,.\nonumber
\end{eqnarray}
In the second step we must show that the above Ricci pinching is
preserved by the ordinary differential equation (\ref{ode}). Let
$\Ric(X_{a,b})$ denote the Ricci tensor of $X_{a,b}$. Assume that
$\lambda_i=-(1-p)\blambda$. We have to show that
\begin{equation*}
\Ric(X_{a,b})_{ii}> p \frac{\scal(X_{a,b})}{n}
=p\Bigl(\bigl(1+2(n-1)a\bigr)\blambda^2
   +\frac{(1+(n-2)b)^2}{1+2(n-1)a}\sigma\Bigr)
\end{equation*}
holds for $b\in (0,1/2]$. We first observe that by (\ref{RicstrRicci})
\begin{eqnarray*}
\Ric(\Rc^2+\Rc^\#)_{ii}
   = \sum_{k\neq i}\Ric_{kk}\Rc_{kiik}
  \ge \sum_{k\neq i} p\blambda\Rc_{kiik}
  = p^2\blambda^2\,.
\end{eqnarray*}
Using formula (\ref{equ}) for $d_{ij}$ and
$\lambda_i=-(1-p)\blambda$ we see that
\begin{eqnarray*}
\Ric(X_{a,b})_{ii}
  &\ge& p^2\blambda^2+ \sum_{j\neq i} d_{ij}\\
  &=& p^2\blambda^2 +2(n-1)ap\blambda^2+ (n-2)b^2 (1-p)^2\blambda^2
      +nb^2 \sigma\\
  & &+\frac{(n-1) \sigma b^2(1-2b)}{(1+2(n-1)a)(1+(n-2)b^2)} \bigl(
     2+2(n-2)b\bigr)\,.
\end{eqnarray*}
By our choice for $b$ and $p$ it is straightforward to check that
$$
 p^2+(n-2)b^2(1-p)^2=p\,.
$$
This shows that in the asserted inequality the
$\blambda^2$-terms cancel each other.
Since $\sigma> 0$ it remains to verify
$$
nb^2+\frac{(n-1) b^2(1-2b)}{(1+2(n-1)a)(1+(n-2)b^2)} \bigl(
2+2(n-2)b\bigr)>p\frac{(1+(n-2)b)^2}{1+2(n-1)a}\,.
$$
The identity
$$
 (1+2(n-1)a)=\frac{1+2(n-1)b+n(n-2)b^2}{1+(n-2)b^2}
$$
shows that this is equivalent to
\begin{eqnarray*}
0&<& n(1+2(n-1)b+n(n-2)b^2)- (n-2) (1+(n-2)b)^2 \\
&&+ (n-1)(1-2b)\bigl( 2+2(n-2)b\bigr)\\
&=&  2n+2n(n-2)b\,.
\end{eqnarray*}
This shows the claim.
\end{proof}

We remark that the above sets remain in fact invariant for all
$b>0$. For $b\to +\infty$ they converge to an invariant set of
Einstein curvature operators.

We will now finish the proof of Theorem~\ref{thm: pinching fam} by showing
 that the cone from
Lemma~\ref{lem: first family} for $b=1/2$ can be joined
by a continuous family of invariant cones with arbitrarily small cones around
the identity.

\begin{lem}\label{lem: seconddef}\label{lem: second family}
Assume $b=1/2$ and put for $s\geq 0$
$$
 a=\frac{1+s}{2}\quad \textrm{ and }  p= 1-\frac{4}{n+2+4s}\,.
$$
Then the set
\begin{eqnarray*}
\li_{a,b}\Bigl(\bigl\{ \Rc\in S^2_B(\so(n))\mid \Rc\ge 0,
   \Ric\ge p(s)\tfrac{\tr(\Ric)}{n}\}\Bigr)
\end{eqnarray*}
is invariant under the vector field corresponding to {\rm
(\ref{ode})}. In fact, it is transverse to the boundary of the set at all
boundary points $\Rc\neq 0$.
\end{lem}

Notice that $\lim_{s\to\infty}\tfrac{1}{a}\li_{a,b}(R)=2(n-1)R_I$. 
Consequently the cones of the lemma converge to $\R^+ I$ for $s\to \infty$.

\begin{proof}
Notice that the formulas in Corollary \ref{cor: Deigen} simplify:
\begin{eqnarray*}
d_{ij}&=&\bigl(\tfrac{1}{4}(n-2)-s\bigr) \lambda_i\lambda_j
+(s+1)( \blambda +\lambda_i)( \blambda +\lambda_j)
+\tfrac{1}{4}(\lambda_i^2+\lambda_j^2)\\
&&-\frac{\sigma ns}{4n+4(n-1)s}
\end{eqnarray*}
and
\begin{eqnarray*}
r_i&=&-\lambda_i^2+ (s+1)\blambda(n-2)
\lambda_i+(s+1)(n-1)\blambda^2 +\frac{\sigma n^2}{4n+4(n-1)s}\,.
\end{eqnarray*}
We first verify that $X_{a,b}$ does preserve the Ricci pinching.
We may suppose that $\lambda_i=-(1-p)\blambda$. We have to show
\begin{eqnarray*}
0&\le& p^2\blambda^2
-(1-p)^2\blambda^2-(s+1)\blambda^2(n-2) (1-p)
+(s+1)(n-1)\blambda^2\\
&&+\frac{\sigma n^2}{4n+4(n-1)s}
- p\Bigl( (n+(n-1)s)\blambda^2+\frac{ n^2}{4n+4(n-1)s}\sigma\Bigr)\,.
\end{eqnarray*}
Because of $\sigma\ge 0$ we can neglect the terms with $\sigma$.
Dividing by $\blambda^2$ gives
\begin{eqnarray*}
 p^2-(1-p)^2+(s+1)+(s+1)p(n-2) - p(n+(n-1)s)
 &=&s(1-p)\,,
\end{eqnarray*}
which is clearly positive. Notice that this calculation is
independent of $p$.
As before we can complete the proof by
showing that $D_{a,b}$ is positive definite. Using
$$
 \sigma
  \le  (n-1)(1-p)^2\blambda^2
   =   \frac{16(n-1)\blambda^2}{(n+2+4s)^2}
$$
we see that
\begin{eqnarray*}
d_{ij}&=&
 \frac{n+2}{4}(\lambda_i +\frac{4\blambda}{n+2})
               (\lambda_j +\frac{4 \blambda}{n+2})
          +s\blambda (\lambda_i+\lambda_j+\frac{8\blambda}{n+2+4s})
     \\
&&  +\frac{1}{4}(\lambda_i^2+\lambda_j^2)
   +\frac{n-2}{n+2}\blambda^2
   +s\frac{n-6+4s}{n+2+4s}\blambda^2
  -\frac{\sigma ns}{4n+4(n-1)s}\\
&\ge& \Bigl(\frac{n-2}{n+2}+s\frac{n-6+4s}{n+2+4s}
-\frac{16(n-1) ns}{(4n+4(n-1)s)(n+2+4s)^2}\Bigr)\blambda^2\\
&>& \bigl(5+s(n-6)+4s^2- 4s\bigr)\frac{\blambda^2}{n+2+4s}\, >\,0
\end{eqnarray*}
where we used $n\ge 3$ in the last two inequalities.
\end{proof}


\section{Constructing a generalized
Pinching Set from a Family of invariant Cones} \label{sec:
pinchset}

We show how to construct from a family of invariant cones a
generalized pinching set, similar to Hamilton's concept in \cite{H2}.
Let us recall that we denoted by $S^2_B(\so(n))$ the space of
curvature operators.

\begin{thm}\label{thm: pinch}\label{thm: bending cones}
Let $C(s)_{s\in [0,1)}\subset S^2_{B}(\so(n))$ be a continuous
family of closed convex $\SO(n)$-inva\-riant cones of full
dimension, such that $C(s)\setminus \{0\}$ is contained in the half
space of curvature operators with positive scalar curvature. Suppose that
for $\Rc \in C(s)\setminus \{0\}$ the vector field
$X(\Rc)=\Rc^2+\Rc^\#$ is contained in the interior of the tangent
cone of $C(s)$ at $\Rc$ for all $s\in (0,1)$. Then for $\eps, h_0 >0$ there
exists a closed convex $\SO(n)$-invariant subset  $F\subset
S^2_{B}(\so(n))$ with the following properties:
\begin{enumerate}
\item $F$ is invariant under the vector field $X$.
\item $C(\eps)\cap \{\Rc \mid \tr(\Rc)\le h_0\}\subset F$.
\item $F\setminus C(s)$ is relatively compact for all $s\in [\eps,1)$.
\end{enumerate}
\end{thm}

We remark that $F$ is $\Or(n)$-invariant if the cones are. We also
note that the analogue of the theorem holds in the vector space of
K\"ahler curvature operators.

\begin{proof}
Let  $F$ denote the minimal closed convex $\SO(n)$-invariant
subset which is invariant under the flow of $X$ and which contains the set
$$
   C(\eps) \cap \{ \Rc \mid \trace(\Rc)\le h_0\}\,.
$$
Notice that $F$ is the intersection of all subsets which satisfy
the above properties. In particular $F$ is well defined and
$F\subset C(\eps)$. We have to prove that for all $s$ the set
$F\setminus C(s)$ is bounded.

Suppose on the contrary that $F\setminus C(s)$ is not bounded for
some $s$. Let $s_0\geq \eps $ denote the infimum among all $s$
with this property. Since the vector field $X$ is transverse to
the boundary of $C(s_0)$ it is clear that for all small $\delta$
the cone $C_\delta(s_0)$ over the convex set
$$
  \{ \Rc \in C(s_0) \mid \tr(\Rc)=1, d(\Rc, \partial C(s_0))\ge \delta \}
$$
is invariant under $X$. For small $\delta$ the cone
$C_\delta(s_0)$ has maximal dimension and $C_0(s_0)=C(s_0)$. We
now choose $\delta_0$ so small that the vector field $X$ is
transverse to the boundary of $C_\delta(s_0)$ for all $\Rc\neq 0$
and for all $\delta\in [0,\delta_0]$.

A simple compactness argument shows that there is some constant
$\eta$ such that for each $R\in C_\delta(s_0)$ the vector field
$X$ has distance at least  $\eta \|R\|^2$ to the boundary of the
tangent cone of $C_\delta(s_0)$ at $\Rc$. We note that $X$ is
locally Lipschitz continuous with a Lipschitz constant that
growths linearly in $\| \Rc\|$. Combining both facts we see that
there is some constant $c>0$ such that the truncated shifted cone
$$
TC_\delta(s_0):=\{ \Rc \mid \Rc +I\in C_\delta(s_0),
  \tr(\Rc)\ge \bh \}
$$
is invariant under the flow of $X$ for all $\delta\in
[0,\delta_0]$.

Consequently for small $\delta>0$ we have
that $C(s_0)\cap \{\Rc \mid \tr(\Rc)=\bh\}$ is contained in the
interior of $TC_\delta(s_0)$.
Since the family $C(s)$ is continuous, we conclude
  \[
C(\bar s)\cap \{\Rc \mid \tr(\Rc)=\bh\}\subset TC_\delta(s_0)
\]
 for some $\eps \le \bar s <s_0$. In the case of $s_0=\eps$ put 
$\bar s=\eps$.
 By the definition of $s_0$ we can choose
$k\in \N$ so large that
\begin{eqnarray*}
   F\cap \{\Rc \mid \tr(\Rc)= k\bh\}\subset
C(\bar s)\cap \{\Rc \mid \tr(\Rc)= k\bh\}\subset k\cdot TC_\delta(s_0)\,.
\end{eqnarray*}
The scaled set $k\cdot TC_\delta(s_0)$ is invariant under the flow of
 $X$ too,
since $X(k\Rc)=k^2X(\Rc)$.
Thus the set
\[
F':=\Bigl(F\cap \{\Rc \mid \tr(\Rc)\le k\bh\}\Bigr)\cup\Bigl( F\cap  k\cdot TC_\delta(s_0)\Bigr)
\]
is convex and invariant under the flow of $X$.
By assumption $F'\subset F$. On the other hand
$
TC_\delta(s_0)\setminus C_{\delta/2}(s_0)$ is bounded.
By the continuity of the family it follows that
$F\setminus C(s)$ is bounded for all $s$ which are sufficiently close to
$s_0$. A contradiction to the choice of $s_0$.
\end{proof}

\section{Proof of the Main Result}

Using Theorem~\ref{thm: pinching fam}, Theorem~\ref{mainthm1} 
 is an immediate consequence
of the following

\begin{thm}\label{thm: from ode to pde}
 Let $C(s)_{s\in [0,1)}\subset S_B^2(\so(n))$  be a pinching
 family of closed convex cones, $n\ge 3$.
Suppose that $(M,g)$  is a compact Riemannian manifold such that
the curvature operator of $M$ at each point is contained in the
interior of $C(0)$.  Then the normalized Ricci flow  evolves $g$
to a constant curvature
 limit metric.
\end{thm}

\begin{proof}
Let $\Rc_p$ denote the curvature operator 
of $(M,g)$ at a point $p\in M$. For all $p\in M$
we have
\[
  \Rc_p \in \{ \Rc \mid \scal\le h_0\}\cap C(\eps)
\]
for a sufficiently small $\eps>0$
and a sufficiently large $h_0$,
since the family of cones is continuous and $M$ is compact.
For this pair $\eps,h_0$ we consider an invariant set $F$ as in
Theorem~\ref{thm: bending cones}. 

By the maximum principle 
the Ricci flow evolves $g$ to metrics $g_t$ whose
curvature operators at each point are contained in $F$. We do also
know that the solution of the Ricci flow exists as long as the
curvature does not tend to infinity. Furthermore it follows from
the maximum principle that the Ricci flow exists only
on a finite time interval $t\in [0,t_0)$. By Shi \cite{Shi} it follows from
the maximum principle applied to the evolution equation for the
$i$-th derivatives of the curvature operator that
\[
\max \| \nabla^i \Rc_t\|^2 \le C_i\max\| \Rc_t\|^{i+2}
\]
for all $t\in [t_0/2,t_0)$.

We now rescale each metric $g_t$ to a metric $\tg_t$ such that
the maximal sectional curvature is equal to $1$. From the above
estimates we get a priori bounds for all
derivatives of the curvature tensor of the metric $\tg_t$ for
$t\in [t_0/2,t_0)$.

Next, we pick a point $p_t\in (M,g_t)$
such that the sectional curvature attains its maximum 
in the ball $B_{\pi}(p_t)$ of radius $\pi$ around $p_t$. 
 We pull the metric via the exponential   map
back to the ball of radius $\pi$ in $T_{p_t}M$.
By choosing a linear isometry $\R^n\to T_{p_t}M$
we identify this ball with the ball $B_{\pi}(0)\subset \R^n$ and
denote by $\bg_t$ the induced metric on $B_{\pi}(0)$.
From the above estimates on the derivatives of the
curvature tensor it is clear that for any sequence
$(t_k)$ in $[0,t_0)$ converging to $t_0$ there is a subsequence of
$(\bg_{t_{k}})$  converging in the $C^\infty$ topology to a limit metric.

Let now $\lambda_j$ denote the scaling factors of these metrics 
$\bg_{t_{j}}$ which by assumption tend to infinity. At each point of $M$
the curvature operator of the limit metric is contained in the set 
$$
  \bigcap \,\,\tfrac{1}{\lambda_j^2} F = \R^+ I\,.
$$
Thus the limit metric on $B_{\pi}(0)$ has pointwise constant
sectional curvature. Since $n\ge 3$, it has constant curvature
one by Schur's theorem.

Since the sequence was arbitrary, the minimal
sectional curvature converges
on a ball of radius $\pi$ around $p_t$ in $(M,g_t)$
to $1$ as well as $t$ tends to $t_0$. 
Notice that this argument works for all $p_t\in B_{\pi}(q_t)$, 
where $q_t$ denotes a point where
 the sectional curvature attains its maximum $1$. 
Therefore the minimal sectional curvature converges on the ball 
of radius $2\pi$ around $q_t$ to $1$ as well.  
By the theorem of Bonnet Myers $\diam(M,\tg_t)\le 3\pi/2$ 
for large $t<t_0$ and
consequently,
also the minimum of the sectional curvature
of $(M,\tg_t)$ tends to $1$ for  $t\to t_0$.

In the case of manifolds one is done since by Klingenberg's
injectivity radius estimate \cite{CE} collapse can not occur. 
Alternatively, one can use the fact that
 $(M,g_t)$ satisfies the  assumption of Huisken's
theorem \cite{Hu} for suitable large $t$.
In the case of orbifolds one has to use additionally 
Proposition~\ref{prop: orbifold} 
from below.
\end{proof}

Let us remark that collapse in the above situation
can also be ruled out by applying Perelman's local
injectivity radius estimate for the Ricci flow \cite{Pe}.

\begin{prop}\label{prop: orbifold} Let $(X,g)$ be a compact 
orbifold with sectional curvature $K$. If $n\geq 3$ and $g$ is strictly
quarter pinched, that is $1/4<K\le 1$,  
then $X$ is the quotient of a Riemannian manifold by a finite 
isometric group action.
\end{prop}

\begin{proof} By replacing $X$ by a cover if necessary we 
may assume that $X$ is not a nontrivial quotient of 
an orbifold by a finite group action. 
We then have to show that $X$ is a manifold.  
Recall that the frame bundle $FX$  of the orbifold $X$, endowed 
with the connection metric of $g$, is a Riemannian manifold.
We consider an $\SO(n)$ orbit $\SO(n)v$ in $FX$. 
Clearly the normal exponential 
map of the orbit $\SO(n)v$ has a focal radius $\ge \pi$.
Similarly to Klingenberg's injectivity radius estimate 
we show below that the normal exponential 
map of the orbit $\SO(n)v$ has injectivity radius $\ge \pi$. 
Since the orbit was arbitrary, this  rules out exceptional orbits 
and hence $X$ is then a manifold. 

From the assumption that $X$ is not a nontrivial quotient 
it follows that the natural map $\pi_1(\SO(n))\rightarrow \pi_1(FX)$
is surjective. This implies that the space $\Omega_{\SO(n)v}FX$ 
of all curves starting and ending in $\SO(n) v$ is connected. 
 The critical levels of the energy functional 
in  $\Omega_{\SO(n)v}FX$ are 
in one to one correspondence to the geodesic loops in the orbifold. 

Suppose on the contrary that the injectivity radius of
the normal exponential map of $\SO(n)v$ is equal to $r<\pi$. 
It is then easy to see that there is a horizontal 
geodesic $c$ of length $2r$ in $\Omega_{\SO(n)v}FX$.  Analogously
 to Klingenberg's long homotopy lemma 
one can 
show that every path $c_s$ in $\Omega_{\SO(n)v}FX$ 
that connects $c=c_0$ with a constant curve $c_1$ 
satisfies $L(c_s)\ge 2\pi$ for some $s$. 
In other words the space of paths of energy $< 2\pi^2$ 
is not connected. 

On the other hand it is straightforward to check 
that the critical points of the energy function with energy 
$\ge 2\pi^2$ have indices at least $n-1\ge 2$.
But then by a standard degenerate Morse theory argument
the loop space $ \Omega_{\SO(n)v}FX$ itself is not connected 
-- a contradiction.
\end{proof}

\begin{rem}
\begin{enumerate}
\item[1.] The main difference between the two-dimensional and
 the higher dimensional case is that in dimension two, Schur's theorem fails.
\item[2.]  Proposition~\ref{prop: orbifold} does not remain valid 
in dimension two either. 
In fact given any positive $\delta<1$, there is a $\delta$ pinched 
two-dimensional orbifold $X$ which is not the quotient of a manifold: 
Consider 
two discs of constant curvature $1$ and with totally geodesic boundary. 
Divide out the cyclic group of order $(p+1)$ from the first disc 
and the cyclic group of order $p$ from the second. 
After scaling the first disc by the factor $\tfrac{p+1}{p}$ 
the two orbifolds can be glued along their common boundary.
By smoothing this example for some large $p$ 
one obtains the claimed result. 
\item[3.] The space of $3$-positive curvature 
operators is not invariant under the ordinary differential equation (\ref{ode})
for $n\ge 4$.
\end{enumerate}
\end{rem}


\hspace*{1em}\\
\begin{footnotesize}
\hspace*{0.3em}{\sc University of M\"unster,
Einsteinstrasse 62, 48149 M\"unster, Germany}\\
\hspace*{0.3em}{\em E-mail addresses: }{\sf cboehm@math.uni-muenster.de}\\
\hspace*{8.5em}{\sf wilking@math.uni-muenster.de}
\end{footnotesize}
\end{document}